\documentclass[twocolumn,10pt]{article}    

\usepackage{cite}
\usepackage{amsmath,amssymb,amsfonts}
\usepackage{algorithmic}
\usepackage{graphicx}
\usepackage{textcomp}
\usepackage{url}

\usepackage{color, xcolor} 
\usepackage{epsfig} 
\usepackage{amsmath, bbm, bm} 
\usepackage{amsfonts, amssymb}
\usepackage{mathrsfs}
\usepackage{anyfontsize}

\usepackage{tikz}
\usetikzlibrary{intersections}
\usepackage{pgfplots}

\usepackage{titlesec}
\usepackage[letterpaper, margin=0.8in]{geometry}

\newcommand {\st}{\text{\footnotesize $(t)$}}
\newcommand {\so}{\text{\footnotesize $(0)$}}
\newcommand {\sta}[1]{\text{\footnotesize $(#1)$}}
\newcommand {\tc}{\tilde{c}}
\newcommand {\vbeta}{\vec{\beta}}

\newcommand {\sS}{{\mathscr{S}}}
\newcommand {\sM}{{\mathscr{M}}}

\newcommand {\cB}{{\mathcal{B}}}

\newcommand {\cI}{{\mathcal{I}}}
\newcommand {\cL}{{\mathcal L}}

\newcommand {\cP}{{\mathcal{P}}}
\newcommand {\cR}{{\mathcal{R}}}
\newcommand {\cS}{{\mathcal{S}}}
\newcommand {\cT}{{\mathcal{T}}}
\newcommand {\cV}{{\mathcal{V}}}

\newcommand {\cY}{{\mathscr{Y}}}

\newcommand {\R} {{\rm I\kern-2.5pt R}}
\newcommand {\C} {{\rm I\kern-5pt C}}

\newtheorem{lemma}{Lemma}

\newtheorem{proposition}{Proposition}
\newtheorem{assumption}{Assumption}

\newtheorem{theorem}{Theorem}
\newtheorem{remark}{Remark}
\newtheorem{example}{Example}
\newtheorem{definition}{Definition}

\newcommand{\beqa}{\begin{eqnarray}}
\newcommand{\eeqa}{\end{eqnarray}}
\newcommand{\beqan}{\begin{eqnarray*}}
\newcommand{\eeqan}{\end{eqnarray*}}
\newcommand{\beq}{\begin{equation}}
\newcommand{\eeq}{\end{equation}}

\newcommand{\bfl}{\begin{flushleft}}
\newcommand{\efl}{\end{flushleft}}
\newcommand{\bgnv}{\begin{verbatim}}
\newcommand{\ednv}{\end{verbatim}}

\newcommand{\myarr}{\begin{array}{lll}}

\newcommand{\bitem}{\begin{itemize}}
\newcommand{\eitem}{\end{itemize}}
\newcommand{\benum}{\begin{enumerate}}
\newcommand{\eenum}{\end{enumerate}}

\newcommand{\orcid}[1]{\href{https://orcid.org/#1}{\textcolor[HTML]{A6CE39}{\aiOrcid}}}

\providecommand{\keywords}[1]
{
  \small	
  \textbf{\textit{Keywords---}} #1
}

\pgfplotsset{compat=1.7}

\title{Epidemic Population Games And Evolutionary Dynamics} 
\author{
    Nuno C. Martins$^a$,
    Jair Cert\'{o}rio$^a$,
    Richard J. La$^a$ 
}
\date{$^a$\small \it{Dept. of ECE and ISR, The University of Maryland, at College Park.}}

\renewcommand\footnotemark{}

\titleformat{\section}
  {\bfseries\fontsize{12}{15}}{\thesection}{1em}{}

\titleformat{\subsection}
  {\normalfont\itshape\fontsize{12}{15}}{\thesubsection}{1em}{}

\let\svthefootnote\thefootnote

\begin{document}


%

%




%

\twocolumn[
\begin{@twocolumnfalse}
\maketitle
\begin{abstract}                          
    \small We propose a system theoretic approach to select and stabilize the endemic equilibrium of an SIRS epidemic model in which the decisions of a population of strategically interacting agents determine the transmission rate. Specifically, the population's agents recurrently revise their choices out of a set of strategies that impact to varying levels the transmission rate. A payoff vector quantifying the incentives provided by a planner for each strategy, after deducting the strategies' intrinsic costs, influences the revision process. An evolutionary dynamics model captures the population's preferences in the revision process by specifying as a function of the payoff vector the rates at which the agents' choices flow toward strategies with higher payoffs. Our main result is a dynamic payoff mechanism that is guaranteed to steer the epidemic variables (via incentives to the population) to the endemic equilibrium with the smallest infectious fraction, subject to cost constraints. We use a Lyapunov function not only to establish convergence but also to obtain an (anytime) upper bound for the peak size of the population's infectious portion.
\end{abstract}
\vspace{5pt}
\hrule
\vspace{6pt}
\end{@twocolumnfalse}
]

\let\thefootnote\relax\footnotetext{
    Corresponding author N.~C.~Martins. This work was supported by 
    AFOSR Grant FA9550-19-1-0315.}
\let\thefootnote\relax\footnotetext{
    Email addresses: \texttt{nmartins@umd.edu} (Nuno C. Martins),
\texttt{certorio@umd.edu} (Jair Certório), \texttt{hyongla@umd.edu}
(Richard J. La).}

\let\thefootnote\svthefootnote

\section{Introduction}

This article has two main tenets: {\bf (i)} We adopt a continuous-time susceptible-infectious-recovered-susceptible
(SIRS)
compartmental epidemic model \cite{Pastor-Satorras2015} in which the aggregate decisions of a population of bounded-rationality agents determine the {\it transmission rate},
which we denote as $\cB\st$ at time $t$. We employ a population game approach in which the agents are nondescript and must choose from a set of available strategies $\{1,\ldots,n\}$. Each strategy will have an effect on $\cB\st$, but the agents' choices are guided by each strategy's payoff, or net reward, resulting from a payoff incentive, or reward, after the intrinsic cost of the strategy is deducted. The collective decision-making of the population follows an evolutionary dynamics model that captures the agents' preferences and assumes that the agents can repeatedly revise their strategies~(see~\S\ref{subsec:EDM}). {\bf (ii)} We formulate and solve a design problem that seeks to steer via payoff incentives the agents' decisions to attain the smallest endemic prevalence of infections, subject to a limit on the long term incentives' cost. The problem envisages dynamic payoff mechanisms whose dynamics can be coupled with the state (epidemic variables) of the SIRS model. We will refer to this coupled system as an {\it epidemic population game}~(see \S\ref{subsec:EPGModel}).

\subsection{Evolutionary Dynamics Model (EDM)}
\label{subsec:EDM}
Each agent follows one strategy at a time, which it can revise repeatedly. A payoff vector $p\st$ in $\mathbb{R}^n$ whose entries quantify the net reward of each strategy influences the revision process, as typically an agent will seek strategies with a higher payoff. We consider $p\st$ as follows:
\begin{equation}
\label{eq:payoffEq}
p\st : =  r\st - c
\end{equation} where $c$ is the vector whose $\ell$-th entry $c_\ell$ is the inherent cost of the $\ell$-th strategy, and $r\st$ is a reward vector meant to incentivize the adoption of safer (costlier) strategies.

Rather than focusing on what each strategy may represent, in our analysis we assume that a vector  $\vbeta$ in $\mathbb{R}_{>0}^n$ is given whose $\ell$-th entry $\vbeta_\ell$ quantifies the effect of strategy $\ell$ towards $\cB\st$ according to:
\begin{equation}
\label{eq:betaaverage}
\cB\st = \vbeta' x\st, \quad t \geq 0
\end{equation} where $x\st$ is the so-called {\it population state} taking values in the standard simplex $\mathbb{X}$ defined below and whose $\ell$-th entry $x_\ell\st$ is the proportion of the population adopting the $\ell$-th strategy at time $t$. 
\begin{equation*}
\mathbb{X}:= \Bigg \{x \in [0,1]^n \ \Big | \ \sum_{i=1}^n x_i =1 \Bigg \}
\end{equation*} 

Following the standard approach in~\cite[Section~4.1.2]{Sandholm2010Population-Game}, the following {\it evolutionary dynamics model} {\bf (EDM)} governs the dynamics of~$x$ in the large-population limit:
\begin{subequations}
\begin{equation}\tag{EDMa} \label{eq:EDM-DEF} \dot {x}\st = \mathcal{V} ( x\st,p\st ), \quad t\geq 0, 
\end{equation}
where the $i$-th component of $\mathcal{V}$ is specified as: 
\begin{multline} \tag{EDMb} \label{eq:EDMfromProtocol} \mathcal{V}_i ( x\st,p\st ) := \underbrace{\sum_{j=1, j \neq i}^{n}  x_j\st \mathcal{T}_{ji} ( x\st,p\st )}_{\text{\footnotesize inflow switching to strategy $i$}}  \\ - \underbrace{\sum_{j=1,j \neq i}^{n} x_i\st \mathcal{T}_{i j} ( x\st,p\st ) }_{\text{\footnotesize outflow switching away from strategy $i$}}
\end{multline}
\end{subequations}

A Lipschitz continuous map $\mathcal{T}: \mathbb{X} \times \mathbb{R}^{n} \rightarrow [ 0,\bar{\mathcal{T}}]^{n \times n}$, with upper bound $\bar{\mathcal{T}} > 0$, is referred to as {\it the revision protocol} and models the agents' strategy revision preferences. In \cite[Part~II]{Sandholm2010Population-Game} and \cite[\S 13.3-13.5]{Sandholm2015Handbook-of-gam} there is a comprehensive discussion on protocols types and the classes of bounded rationality rules they model. In \cite[\S{IV}]{Park2019From-Population}, the authors substantiate using (EDM) as a deterministic approximation for the case when a dynamical payoff mechanism generates $p$ from $x$, as will be the case~here. 

Below, we define a widely-used class of protocols, which we will repeatedly invoke to illustrate key concepts and employ in examples throughout the article.

\begin{definition} \label{def:IPC} Any protocol is said to be of the {\it impartial pairwise comparison} (IPC) type~\cite{Sandholm2010Pairwise-compar} if there is a map ${\phi:\mathbb{R}_{\geq 0} \rightarrow [0,\bar{\cT}]^n}$, whose components satisfy $\phi_j(0)=0$ and $\phi_j(\nu)>0$ for $\nu>0$, such that $\cT$ can be recast as:
\begin{equation}
\label{eq:defIPC}
\mathcal{T}_{ij}(x,p) \underset{\text{\tiny IPC}}{=} \phi_j([\tilde{p}_{ij}]_+),
\end{equation} where $\tilde{p}_{ij}:=p_j-p_i$. The well-known Smith's protocol originally proposed to model the commuters' preferences in traffic assignment problems~\cite{Smith1984The-stability-o} can be specified by ${\phi_j^{\text{\tiny Smith}}([\tilde{p}_{ij}]_+) := \min \{ \lambda [\tilde{p}_{ij}]_+, \bar{\mathcal{T}} \} }$. Hence, according to Smith's protocol, the rate of switching from the $i$-th strategy to the $j$-th strategy is proportional to the positive part of the payoff difference $\tilde{p}_{ij}$, up to the upper bound $\bar{\mathcal{T}}$.
\end{definition}

\subsection{Epidemic Population Game (EPG)}
\label{subsec:EPGModel}

We assume that, within the time interval of interest, the population's size is $N\st=e^{gt}N\sta{0}$, where $g$ is a constant representing the difference between the birth and death rates. Here, $N\sta{0}$ is large and $N\st$ approximates  the population's cardinality at time $t\geq0$.
Below is our definition of an {\it epidemic population game} (EPG):
\begin{align} \tag{EPGa}
\dot{I}\st &= \big ( \cB \st (1-I\st-R\st) -\sigma) I\st \\ \tag{EPGb}
\dot{R}\st &= \gamma I\st - \omega R\st \\ \tag{EPGc}
\dot{q}\st & = G (I\st,R\st,x\st,q\st) \\ \tag{EPGd}
r\st & = H(I\st,R\st,x\st,q\st)
\end{align} where $I\st$, $R\st$ and $S\st:=(1-I\st-R\st)$ take values in $[0,1]$ and represent the proportions of the population which are infectious, have recovered and are susceptible to infection at time $t$, respectively. Specifically, these variables are the numbers of infectious, recovered and susceptible individuals at time $t$ divided by $N\st$. Hence, (EPGa,b) is a normalized SIRS model specified by the constants $\sigma:=g+\bar{\sigma}$, $\omega:=g+\bar{\omega}$ and $\gamma<\bar{\sigma}$, which we assume are all positive. We also assume that newborns are susceptible and the disease death rate associated with the epidemic is zero (death rate is independent from the epidemic). {This is a reasonable assumption when the number of deaths caused by the disease is negligible relative to that from all other causes.} In this case, $\gamma$ is the recovery rate, $\bar{\sigma} - \gamma$ is the death rate, and $\bar{\omega}$ is the rate at which recovered individuals become susceptible (due to waning immunity) or expire. Our time unit is \underline{one day}, and $\bar{\sigma}^{-1}$ is the mean infectious period (in days) for an affected individual. We seek to design the dynamic payoff mechanism (EPGc,d), where $r\st$ and $q\st$ take values respectively in $\mathbb{R}^n$ and $\mathbb{R}^m$, with $m \geq 1$, and $G$ and $H$ are Lipschitz continuous.

\subsection{Problem formulation and paper structure}
\label{sec:Problem}

The strategies' inherent costs decrease for higher transmission rates, and we order the entries of $\vbeta$ and $c$~as:
\begin{equation*}
 \vbeta_i < \vbeta_{i+1} \text{ and } c_i > c_{i+1}, \quad 1 \leq i \leq n-1
\end{equation*} We consider that $\vbeta_1 > \sigma$, i.e., a transmission rate less than or equal to $\sigma$ would be unfeasible or too onerous. 

{\bf Convention:} Henceforth, $c$ and $\vbeta$ satisfying the conditions above are assumed given and fixed. Hence, we can simplify our notation by omitting $c$ and $\vbeta$ from this point onward. We will also use $\tilde{c}$ defined below:
\begin{equation*}
    \tilde{c}_i=c_i-c_n, \quad 1 \leq i \leq n
\end{equation*}

\begin{definition}  \label{def:betastar}  Given a cost budget $c^*$ in $(0,\tilde{c}_1)$, we determine the optimal endemic transmission rate $\beta^*$ as:
 \begin{equation} \label{eq:betadef}
\beta^* : = \min \big \{ \vbeta' x \ | \ \tilde{c}'x \leq c^*, \ x \in \mathbb{X} \big \} 
\end{equation}
\end{definition}

\noindent {\bf Main Problem:} We seek to obtain Lipschitz continuous $G$ and $H$ for which the following hold for any  $I\so$ in $(0, 1]$, $R\so$ in $[0,1-I\so]$, $x\so$ in $\mathbb{X}$, and $q\so$ in~$\mathbb{R}^m$:
\begin{align} \tag{P1}
\lim_{t \rightarrow \infty} (I,R,\cB)\st &= (I^*,R^*,\beta^*)  \\ \tag{P2}  \limsup_{t \rightarrow \infty} x\st ' r\st &\leq c^*  ,
\end{align} where, from Picard's Theorem, $\{(I,R,x,q)\st \ | \ t\geq 0\}$ is the unique solution of the initial value problem for the closed loop system formed by~(EDM) and~(EPG). Here, the nontrivial endemic equilibrium for (EPGa,b)~is:
\begin{equation*}
I^* : = \eta(1 - \tfrac{\sigma}{\beta^*}), \ R^* : = (1-\eta) (1 - \tfrac{\sigma}{\beta^*}), \quad \eta:= \tfrac{\omega}{\omega+\gamma}
\end{equation*} 
We will seek $G$ and $H$ for which a \underline{Lyapunov function} \underline{for the overall system} exists. We will do so not only to establish (P1) but, crucially, also to leverage the Lyapunov function to obtain \underline{anytime upper bounds} for $I\st$. This is relevant because, as has been pointed out in studies~\cite{Godara2021A-control-theor,Sontag2021An-explicit-for} employing $\cB\st$ as a control variable, $I\st$ tends to significantly \underline{overshoot} its endemic equilibrium $I^*$ when $I\sta{0}<I^*$, unless the control policy prevents it.

 \begin{remark} \label{rem:naive} Had we \underline{not required a Lyapunov function}, we could have used $r\st=\tilde{c}+\mu(\check{x}-x\st)$, which is a simple potential game~\cite{Monderer1996Potential-games} decoupled from the epidemic variables. Here, $\mu$ could have been any positive constant and $\check{x}$ any solution of~(\ref{eq:betadef}). Such a na\"{i}ve approach would have guaranteed, for a large class of protocols~\cite[\S13.6.5]{Sandholm2015Handbook-of-gam}, that $\lim_{t \rightarrow \infty}\cB\st=\vbeta'\check{x}=\beta^*$. Even if~(P1)-(P2) could have ultimately held (we did not prove it), this approach would have imparted no useful anytime bounds for~$I\st$.
 \end{remark}

\begin{remark} We interpret $\int_t^{t+T} r\sta{\tau}'x\sta{\tau}d\tau$ as the normalized cost of using $r\st$ in the interval $[t,t+T]$, for $T>0$. Hence, (P2) would guarantee that the long-term normalized time-averaged cost a social planner would have to bear for employing $G$ and $H$ would not exceed $c^*$. Moreover, since $I^*$ is an increasing function of $\beta^*$, a solution satisfying (P1) would guarantee that the size of the population's endemic infectious portion would be the lowest, subject to~(P2).
\end{remark}

\noindent {\bf Paper structure:} In~\S\ref{sec:MotAndComp} we motivate our paradigm and compare it with previous work. In~\S\ref{sec:solution} we describe a choice for $G$ and $H$, introduce a suitable Lyapunov function, and state Theorem~1 asserting that our choice (of $G$ and $H$) solves our Main~Problem. We show in~\S\ref{sec:AnytimeBounds} how to use the Lyapunov function to construct anytime bounds, which we also validate numerically via simulation. The article ends with brief conclusions in~\S\ref{sec:conclusions}, and in Appendix~\ref{sec:proof} we rigorously prove Theorem~1.

\section{Motivation and comparison to prior work}
\label{sec:MotAndComp}

Modern theoretical epidemiology can be
traced back to the studies by 
En'ko~\cite{Enko1889}, Hamer~\cite{Hamer1906}, 
Ross~\cite{Ross1908}, and Kermack and McKendrick
\cite{Kermack1927A-contribution-, 
Kermack1932Contributions-t, cKend}. In particular, 
Kermack and McKendrick
used a deterministic model for modeling the transmissions 
in a closed population, which is now known as the 
susceptible-infected-recovered (SIR) model, and demonstrated 
the existence of a critical threshold density of susceptible 
individuals for the occurrence of a major epidemic.
Since then, many related compartmental models have 
been introduced with additional states, e.g.,  
deceased (D), exposed (E), 
maternally-derived (M), vaccinated 
(V), and include 
SEIR/S, SIRD, SIRV, SIS, and MSIR, in addition 
to the SIRS model adopted for our study. 
A comprehensive survey can be found in 
a manuscript by Anderson and May
\cite{Anderson1991}. 

An important aspect of epidemic processes is
human behavior and the strategic interactions 
among individuals, which determine their decisions 
over time in response to their (perceived) payoffs
and in turn shape the course of epidemic process. 
Game theory provides a natural framework and tools for 
studying such strategic interactions, and 
several recent studies adopted an {\em 
evolutionary} or {\em population game} framework
\cite{Amaral2021, Arefin2020, Bauch2003, Bauch2004, 
Bauch2005, Chang2019, Cho2020, d'Onofrio2011,
Kabir2019b, Kabir2020, Kuga2018, Wang2020}.
We refer an interested reader to 
\cite{Chang2020} and references therein for
a comprehensive survey of earlier studies and
a more detailed summary of studies reported
in \cite{Bauch2004, Bauch2005, d'Onofrio2011}.

For example, 
Amaral et al.~\cite{Amaral2021} studied the effects
of perceived risks, i.e., individual cost from
infection, when individuals can choose to 
voluntarily quarantine or continue their normal life. 
They showed that increased perceived risks result in
multiple infection peaks due to strategic interactions.
Kabir and Tanimoto~\cite{Kabir2020} considered
a similar setting and showed that naturally
acquired shield immunity is unlikely to be 
effective in suppressing an epidemic without 
additional social measures with low costs for
individuals. 

In another line of related research, 
di Lauro et al.~\cite{diLauro2021}
and Sontag~\cite{Sontag2021An-explicit-for}
studied the problem of identifying when 
non-pharmaceutical interventions (NPIs), such 
as quarantine and lockdowns, should be
put in place to minimize the peak infections;
\cite{diLauro2021} studied the optimal timing
for one-shot intervention, whereas 
\cite{Sontag2021An-explicit-for} considered
a fixed number of complete lockdowns. 
Al-Radhawi et al. 
\cite{Al-Radhawi2021} modeled media coverage, 
public health measures and other NPIs during 
a prolonged epidemic 
as feedback effects and examined the problem 
of tuning NPIs to regulate infection rates 
as an adaptive control problem. Using 
a singular-perturbation approach, the authors
investigated the stability of disease-free
and endemic steady states. 
{Godara et al. 
\cite{Godara2021A-control-theor} considered the problem 
of controlling the infection rate to minimize the total
cost till herd immunity is achieved in an SIR model. 
They formulated it as an optimal control problem subject 
to a constraint on the fraction of infectious
population.}

{Although we do not consider epidemic processes on 
general networks in this paper, their dynamics on 
networks have been studied extensively
(see \cite{Pare2020, Pastor-Satorras2015} for a review
of the literature), including time-varying
networks~\cite{Pare-TCNS18} and the influence of 
network properties on epidemic processes
 \cite{La2018, La2019}.
Recently, the topic of mitigating disease 
or infection spread in a network has enjoyed 
much attention. In particular, researchers
investigated optimal strategies using
vaccines/immunization (prevention) 
\cite{Cohen-PRL03, Preciado-CDC13}, antidotes 
or curing rates (recovery) \cite{Borgs-RSA10,  
mai2018distributed, Ottaviano-JCN18} 
or a combination of
both preventive and recovery measures
\cite{Nowzari-TCNS17, Preciado-TCNS14}. 
For example, \cite{Preciado-CDC13}
studied the problem of partial vaccination
via investments at each individual 
to reduce the infection rates, with the
aim of maximizing the exponential decay 
rate to control the spread of an epidemic. 
Similarly, \cite{mai2018distributed} 
examined the problem of determining the
optimal curing rates for distributed
agents under different formulations.}

{Our study advances the state-of-the-art 
in several directions: unlike the studies that 
aim to suppress epidemic spread
\cite{Cohen-PRL03, Preciado-CDC13, Borgs-RSA10,
mai2018distributed, Ottaviano-JCN18, 
Nowzari-TCNS17, Preciado-TCNS14}, our goal is to
design policies for minimizing the endemic transmission
rate subject to a constraint on the
long-term average cost a planner bears. Moreover, 
even though \cite{diLauro2021, Godara2021A-control-theor,
Sontag2021An-explicit-for} investigated a related
problem of managing infection rates during epidemics,  
these studies did not consider strategic interactions 
among many agents of bounded rationality, 
which can revise their strategies over time, leading
to more complex dynamics. 
Finally, to the best of our knowledge, 
our study is the first to provide a methodology
for designing policies that can guarantee (a)
provable convergence to an equilibrium set
(see Theorem~\ref{thm:MainTheorem} and Remark 
\ref{rem:thm1Universality}) and (b) fulfill
an {\em anytime} bound on $I(t)$ (see eq. \eqref{eq:AnytimeBoundnis2} and 
Remark~\ref{rem:AnytimeConstraint}). 
As we will discuss in detail, it is notable 
that these results hold under {\em any}
revision protocol $\mathcal{T}$ that satisfies
some assumptions stated in the subsequent section
{\em without the need to know the exact protocol}. 
}


\section{A Solution to Main Problem}
\label{sec:solution}

In this section, we will propose a choice of $G$ and $H$ for (EPG) and in \S\ref{subsec:MainResult} we will state Theorem~1 (our \underline{main result}), which guarantees that our choice solves the Main~Problem, as stated in~\S\ref{sec:Problem}. 

\subsection{Cases~I and~II, and determining $x^*$}
Before we proceed, we will introduce a  definition, an assumption and a related remark.

\begin{definition} \label{def:cases} {\bf (Cases~I~and~II)} Given $c^*$ in $(0,\tc_1)$, one of two cases holds: {\bf Case~I} is defined by when $\tc_{i^*+1} < c^* < \tc_{i^*}$ for some positive $i^* \leq n-1$. {\bf Case~II} occurs when $n \geq 3$ and $c^* = \tc_{i^*}$ for some positive $i^*\leq n-1$.
\end{definition}

\begin{assumption} \label{assm:AboutBeta}
The following must hold when~$n \geq 3$:
\begin{equation}
\label{eq:BetaAndCIneq}
\frac{c_i-c_{i+1}}{\vbeta_{i+1}-\vbeta_i} > \frac{c_{i+1}-c_{i+2}}{\vbeta_{i+2}-\vbeta_{i+1}}, \quad 1 \leq i \leq n-2
\end{equation}
\end{assumption}

According to~(\ref{eq:BetaAndCIneq}), we assume that as the transmission rate decreases it becomes costlier to reduce it further.
\begin{remark} \label{rem:KKTxstar} Subject to Assumption~\ref{assm:AboutBeta}, it follows from a straightforward application of Karush–Kuhn–Tucker conditions that, for any given $c^*$ in $(0,\tc_1)$,~(\ref{eq:betadef}) has a unique solution we denote as:
\begin{equation}
    x^*: = \arg \min \big \{ \vbeta' x \ | \ \tc'x \leq c^*, \ x \in \mathbb{X} \big \} 
\end{equation} For Case~I, with $\tc_{i^*+1} < c^* < \tc_{i^*}$, it results that $x^*_{i^*} = (c^*-\tc_{i^*+1})/(\tc_{i^*}-\tc_{i^*+1})$, $x^*_{i^*+1}=1-x^*_{i^*}$ and the other entries of $x^*$ are zero, while for Case~II, with $c^* = \tc_{i^*}$, we get that $x^*_{i^*}=1$ and the other entries of $x^*$ are zero. We also immediately conclude from $\beta^*=\vbeta'x^*$ that $\vbeta_{i^*} < \beta^* < \vbeta_{i^*+1}$ for Case~I and $\beta^* = \vbeta_{i^*}$ for Case~II.
\end{remark}

\subsection{A choice of $G$ and $H$}
\label{subsec:definitions}
We start by defining $(\hat{I},\hat{R})$ below, which can be interpreted as {\it ``reference"} epidemic variables determined by the population state $x$ via $\cB=x'\vbeta$:
\begin{equation} \label{eq:referenceRandI}
\hat{I} : =  \eta \Big( 1 - \frac{\sigma}{\cB} \Big), \ 
\hat{R} : =  (1-\eta) \Big( 1 - \frac{\sigma}{\cB} \Big)
\end{equation}

A soon to be described Lyapunov function motivated the following choice of $G$ and $H$:
\begin{subequations}
\label{eq:GHchoice}
\begin{align} \nonumber
G(I,R,x,q) : =& \ (\hat{I} - I) + \eta( \ln I - \ln \hat{I} )  + \upsilon^2 (\beta^* -\cB) \\&+\tfrac{\cB}{\gamma}(R-\hat{R})(1-\eta-R) \label{eq:GDef} \\
H(I,R,x,q) :=& q \vbeta + r^* \label{eq:HDef} 
\end{align} 
where $\upsilon>0$ and $\rho^* > 0$ (see~\S\ref{subsubsec:RulesForRho}) are design parameters, and $r^*$ is the following stationary reward vector:
\begin{equation*} 
    r^*_i := \begin{cases} \tc_i-\rho^* & \text{if $x^*_i=0$} \\ \tc_i & \text{otherwise} \end{cases}, \quad 1\leq i \leq n 
\end{equation*}
\end{subequations}
In \S\ref{subsubsec:RulesForRho}, we describe the rules for selecting a valid $\rho^*$.

\begin{definition}{\bf (Design parameters)} We refer to ${\upsilon >0}$, $\rho^* > 0$ and $c^*$ in $(0,\tc_1)$ as { design parameters}. Specifically, $\upsilon^2$ appears in $G$ and $H$, $\rho^*$ is present in $r^*$, and $c^*$ determines $\beta^*$ and $x^*$.
\end{definition}

In Appendix~\ref{sec:proof}, our proof for the upcoming Theorem~1 will use the fact that, for the $r^*$ chosen, $x^*$ will be the only element $x$ of $\mathbb{X}$ that simultaneously satisfies $\vbeta'x = \beta^*$, and also maximizes $x'(r^*-c)$, which is equivalent to it being the best response to the equilibrium payoff $r^*-c$. 

Notice that, when the epidemic is beginning or is effectively contained, (\ref{eq:GDef}) can be approximated simply as:
 \begin{equation*} 
     G(I,R,x,q) \underset{\text{\tiny $R,I \ll 1$}}{\approx}  \eta ( \ln I - \ln \hat{I})  + \upsilon^2 (\beta^* -\cB)
 \end{equation*}
 According to~(EPGc), $G$ will govern the dynamics of $q\st$, which will indirectly regulate $\cB\st$ via the payoff $p\st$ in~(\ref{eq:payoffEq}) and $H$ in~(EPGd). Specifically, if the population's agents adhere to a protocol, such as IPC, that prioritizes strategies with higher payoffs, then lowering $q\st$ would have the effect of decreasing more the payoff of riskier strategies and hence steering the population towards safer strategies that lower $\cB\st$ and $\hat{I}\st$. On the contrary, increasing $q \st$ would incentivize higher~$\cB\st$.

\subsubsection{Rules for selecting a valid $\rho^*$} \label{subsubsec:RulesForRho}

We start by observing that if $n=2$ then, since $c^*$ is in $(\tilde{c}_2,\tilde{c}_1)$, we have Case~I and from Remark~\ref{rem:KKTxstar} we can further conclude that $\rho^*$ is not present in $r^*$.

To select a valid $\rho^* > 0$ for $n \geq 3$, proceed as follows: {\bf (i)}~For Case~I, choose any $\rho^* >0$. {\bf (ii)}~For Case~II, select any $\rho^* \geq \max \{\vbeta_n - \beta^*,\beta^* - \vbeta_1 \}$.


\subsection{Nash stationarity and $\delta$-passivity assumption}
\label{subsec:DeltaPass}

\begin{definition} {\bf (Nash Stationarity)}
A protocol $\cT$ is ``Nash~stationary" if the following holds for all $p$ in $\mathbb{R}^n$:
\begin{equation}\tag{NS}
\cV(x,p)=0 \quad \Leftrightarrow \quad x \in \sM(p),
\end{equation} where $\sM : \mathbb{R}^n \rightarrow 2^{\mathbb{X}}$ is the following best response map\footnote{A best response in our context may be interpreted in the mass-action sense~\cite{Jr.1951Non-Cooperative} also discussed in~\cite{Weibull1995The-mass-action}.}:
\begin{equation*}
\sM(p) : = \arg \max_{x \in \mathbb{X} } \ p'x, \quad p \in \mathbb{R}^n. \end{equation*}
\end{definition}

Therefore, for a protocol satisfying (NS), $x$ is an equilibrium of~(EDM) if and only if $x$ is a best response to~$p$. Any IPC protocol (see~Definition~\ref{def:IPC}) satisfies~(NS), as do other large classes of protocols~(see~\cite[\S13.5.3]{Sandholm2015Handbook-of-gam}).


Our analysis of the long-term evolution of $(\cI,\cR)\st$ and $(x,p)\st$ will leverage the following assumption stemming from the $\delta$-passivity concept originally proposed in~\cite{Fox2013Population-Game} and later generalized in~\cite{Park2019From-Population,Arcak2020Dissipativity-T,Kara2021Pairwise-Compar-a}.

\begin{assumption}\label{assm:deltaD} 
There exist a differentiable function ${\cS: \mathbb{X}\times \mathbb{R}^n \rightarrow \mathbb{R}_{\ge 0}}$ and a Lipschitz continuous function ${\cP: \mathbb{X}\times \mathbb{R}^n \rightarrow \mathbb{R}_{\ge 0}}$ that satisfy the following inequality for all $x$, $p$ and $u$ in $\mathbb{X}$, $\mathbb{R}^n$ and $\mathbb{R}^n$, respectively:
\begin{subequations}
\label{eq:delta-passivityWconditions}
\begin{equation}\label{delta-passivity}
\frac{\partial \cS(x,p)}{\partial x}\cV(x,p)+\frac{\partial \cS(x,p)}{\partial p}u \le -\cP(x,p)+u' \cV(x,p) 
\end{equation}
where $\cS$ and $\cP$ must also satisfy the equivalences below:
\begin{align}
\label{informative}
\cS(x,p)=0 \quad &\Leftrightarrow \quad \cV(x,p)=0. \\
\label{negdef}
\cP(x,p)=0 \quad &\Leftrightarrow \quad \cV(x,p)=0 
\end{align}

In addition, the following inequality (not required in standard $\delta$-passivity) must hold:
\begin{equation}
\label{homogeneity}
\cP(x,\alpha p) \geq \cP(x,p), \quad \alpha \geq 1, x \in \mathbb{X}, p \in \mathbb{R}^n
\end{equation}
\end{subequations}
\end{assumption}

Later on, in Appendix~\ref{sec:proof},~(\ref{homogeneity}) will be useful to cope with the lack of an a-priori bound for $|p\st|$. 

Based on the Lyapunov functions in~\cite{Hofbauer2009Stable-games-an}, the authors of~\cite{Fox2013Population-Game} determined, for main classes of protocols, explicit expressions for $\cS$ and $\cP$, of which the following is an important example. Explicit constructions for $\cS$ and $\cP$ for 
a generalization of IPC protocols can be found in~\cite{Kara2021Pairwise-Compar-a}.

\begin{subequations}
\begin{remark} \label{rem:StructureForIPC} For any IPC protocol~(\ref{eq:defIPC}) with non-decreasing $\{\phi_1, \ldots,\phi_n\}$, the following satisfy~(\ref{eq:delta-passivityWconditions}):
\begin{align}
\label{eq:cSIPCFormula}
\cS^{\text{\tiny IPC}}(x,p) &:=  \sum_{i=1,j=1}^{n,n} x_i  \int_0^{[\tilde{p}_{ij}]_+}\phi_j(\nu)d\nu \\
\cP^{\text{\tiny IPC}}(x,p) &:=  -\sum_{i=1,j=1}^{n,n} \cV_i^{\text{\tiny IPC}}(x,p)  \int_0^{[\tilde{p}_{ij}]_+}\phi_j(\nu)d\nu
\end{align} where $\cV_i^{\text{\tiny IPC}}$ is obtained by substituting the IPC protocol into~(EDMb). The analysis in~\cite[below the expression for $\dot{\Psi}(x)$ on p.1691]{Hofbauer2009Stable-games-an}
can be
used here to show that $\cP^{\text{\tiny IPC}}$ is non-negative and satisfies~(\ref{delta-passivity}) and~(\ref{negdef}). A small modification of the same argument shows that~(\ref{homogeneity}) holds when $\{\phi_1, \ldots,\phi_n\}$ are non-decreasing.
\end{remark} 
\end{subequations}

\subsection{A Lyapunov function and its properties}
\label{subsec:Lyapunov}

We start with the following reparameterization:
\begin{subequations}
\label{eq:param}
\begin{align} \label{eq:param-a}
&\cI\st := \cB\st I \st & &\cR\st := \cB \st R \st \\
&\hat{\cI}\st : = \cB\st \hat{I}\st & &\hat{\cR}\st : = \cB\st \hat{R}\st 
\end{align}
\end{subequations} 
where $\hat{I}(t)$ and $\hat{R}(t)$ are defined 
in \eqref{eq:referenceRandI}.
Using~(\ref{eq:GHchoice}) and~(\ref{eq:param}), we  rewrite (EPG) as follows:
\begin{subequations}
\label{eq:NewEPG}
\begin{align} \label{eq:NewEPGa}
 \dot{\cI}\st =& \ \cI\st \big ((\hat{\cI}\st-\cI\st) + (\hat{\cR}\st-\cR\st)  \big ) + I\st\dot{\cB}\st \\ \label{eq:NewEPGb}
\dot{\cR}\st =&  \ \omega(\hat{\cR}\st-\cR\st) - \gamma (\hat{\cI}\st- \cI\st ) +R\st \dot{\cB}\st \\ 
\dot{q}\st  =& \ G(I\st,R\st,x\st,q\st)  \label{eq:NewEPGc} \\
r\st  =& \ H(I\st,R\st,x\st,q\st) = q\st \vbeta + r^*, \label{eq:NewEPGd}
\end{align} which after substitution into~(\ref{eq:payoffEq}) gives:
\begin{equation}
p\st = q\st \vbeta + r^o, \quad r^o:= r^*-c \label{eq:NewEPGe}
\end{equation}
\end{subequations}

Here, we modified the representation of the SIRS model in~\cite[between~(3) and~(4)]{ORegan2010Lyapunov-functi} to obtain~(\ref{eq:NewEPGa})-(\ref{eq:NewEPGb}).

We now proceed to define the Lyapunov function:
\begin{align}
\label{eq:defcL}
\cL(\cY) : =&  \cS(x,p) + \sS(\cI,\cR,\cB), \quad \cY \in \mathbb{Y}, \\ \nonumber
\cY:=& (\cI,\cR,x,q), \quad \cB:=\vbeta'x
\end{align} where $\sS$ is defined below, $\cS$ satisfies~(\ref{eq:delta-passivityWconditions}), and $\cY$ taking values in $\mathbb{Y}$ defined below is the state of the complete system comprising~(EDM) and~(\ref{eq:NewEPG}) :
\begin{multline}
\mathbb{Y}:= \big \{(\cI,\cR,x,q) \ | \\ \vbeta_1 \leq \cB \leq \vbeta_n, 0 < \cI \leq \cB, 0\leq \cR \leq \cB - \cI,\\ x \in \mathbb{X}, q \in \mathbb{R}^n \big \}.
\end{multline} Here,~$\sS$ defined below is a modification of the elegant Lyapunov function in~\cite{ORegan2010Lyapunov-functi}:
\begin{equation*} 
\sS(\cI,\cR,\cB) : = (\cI-\hat{\cI}) + \hat{\cI} \ln \tfrac{\hat{\cI}}{\cI} +\tfrac{1}{2\gamma}(\cR-\hat{\cR})^2 + \tfrac{\upsilon^2}{2}(\cB - \beta^*)^2
\end{equation*}
\noindent {\bf Notation convention:} We note that $\sS$ depends on the given design parameters $c^*$, $\upsilon$ and $\rho^*$. But, to simplify our notation we decided not to indicate this dependence.

\begin{remark} \label{rem:PropOfsS} Notice that $\sS$ is convex and nonnegative, and $\sS(\cI,\cR,\cB) = 0$ if and only if $(\cI,\cR,\cB) = (\cI^*,\cR^*,\beta^*)$, where $\cI^* : = \beta^* I^*$ and $\cR^* : = \beta^* R^*$.
\end{remark}

Consequently,~(\ref{eq:delta-passivityWconditions}) and~(NS) imply that $\cL$ satisfies:
\begin{multline} \label{eq:LeqZero}
\cL(\cY) = 0 \Leftrightarrow \\ \Big ((\cI,\cR,\cB) = (\cI^*,\cR^*,\beta^*) \text{ and } x \in \sM(q\vbeta + r^o) \Big )
\end{multline} 

After taking derivatives, and using~(\ref{eq:delta-passivityWconditions}) and~(\ref{eq:NewEPG}), we get:
\begin{align} 
\tfrac{d}{dt} \cL(\cY\st) \leq &- \cP(x\st,p\st) \nonumber  \\  &-(\cI\st-\hat{\cI}\st)^2 -\frac{\omega}{\gamma}(\cR\st-\hat{\cR}\st)^2  \label{eq:DerLIneq}
\end{align}

\subsection{A stability concept and main result}
\label{subsec:MainResult}

Despite the technical issues discussed in \S\ref{subsec:TechnicalIssues}, we will be able to leverage $\cL$ and~(\ref{eq:DerLIneq}) to establish global asymptotic stability of an equilibrium set in the following sense.
\begin{definition} \label{def:GAS}  A set $\mathbb{E} \subset \mathbb{Y}$ is said to be globally asymptotically stable {\bf (GAS)} if it satisfies the conditions (GASa)-(GASc) for the Lyapunov function $\cL$ in~(\ref{eq:defcL}) and the feedback system formed by~(EDM) and~(\ref{eq:NewEPG}):
\begin{itemize} 
    \item[] {\bf (GASa)} It holds that $\cY \in \mathbb{E}\Leftrightarrow \cL(\cY) = 0$.
    \item[] {\bf (GASb)} Given any $\cY\sta{0}$ in $\mathbb{Y}$, $\{ \cY\st \ | \ t\geq 0\}$ has at least one accumulation point in~$\mathbb{E}$.
    \item[] {\bf (GASc)} Given any $\cY\sta{0}$ in $\mathbb{Y}$, \underline{all} accumulation points of $\{ \cY\st \ | \ t\geq 0\}$ are in $\mathbb{E}$.
\end{itemize}
\end{definition}



The following theorem is our main result, in which statement we use Cases~I and~II, as specified in Definition~\ref{def:cases}. Appendix~\ref{sec:proof} provides a rigorous proof of the theorem and  \S\ref{subsec:TechnicalIssues} discusses relevant technical issues.

\begin{theorem} \label{thm:MainTheorem} {\bf (Main result)} Let the protocol defining~(EDM) and the design parameters $\upsilon>0$, $\rho^* > 0$ (valid according to \S\ref{subsubsec:RulesForRho}) and $c^*$ in $(0,\tc_1)$ be given. If (NS) and Assumptions~\ref{assm:AboutBeta}-\ref{assm:deltaD} hold, then the set $\mathbb{E}^*$ defined below is GAS:
\begin{align}    \nonumber
\mathbb{E}^* : =& (\cI^*,\cR^*,x^*) \times \mathfrak{Q}^*, \quad
\mathfrak{Q}^*: = \begin{cases} \{0\} &\text{\small(Case~I)} \\ [-\zeta_2^*,\zeta_1^*] &\text{\small(Case~II)}
\end{cases}
\end{align} where $\zeta_1^*:= \rho^*(\vbeta_n-\beta^*)^{-1}$, and $\zeta_2^*:=  \rho^*(\beta^*-\vbeta_1)^{-1}$.
\end{theorem}
According to Theorem~1, for any $\cY\sta{0}$ in $\mathbb{Y}$, $\cY\st$ will converge to $\mathbb{E}^*$. Consequently, given that $\vbeta_n \geq \cB\st \geq \vbeta_1 >0$, from~(\ref{eq:param}) we conclude that the following~holds:
\begin{equation}
  (I,R,x,q)\st \underset{t \rightarrow \infty}{\longrightarrow} (I^*,R^*,x^*) \times \mathfrak{Q}^*
\end{equation} implying~(P1) as defined in our Main~Problem. This will also imply~(P2) for (Case~I). For (Case~II), we infer { $\lim_{t \rightarrow \infty} r\st'x\st \leq c^* + \zeta_1^*\beta^* $}. Since (Case~I) holds for any $c^*$ in $(0,\tilde{c}_1)$, except for a finite set of values, it is unimportant that (P2) is not guaranteed for (Case~II).
\begin{remark} \label{rem:thm1Universality} {\bf (Universality of Theorem~1)} According to Theorem~1, the payoff mechanism~(\ref{eq:NewEPGc})-(\ref{eq:NewEPGe}) will guarantee that $\mathbb{E}^*$ is GAS for \underline{any} protocol satisfying~(NS) and Assumption~\ref{assm:deltaD}. This is relevant because even when the exact protocol is unknown, one may still be able to conclude from its structure that it satisfies~(NS) and Assumption~\ref{assm:deltaD}. The IPC protocol class is a case in point as (NS) and Assumption~\ref{assm:deltaD} hold for \underline{any} $\bar{\cT}>0$ and \underline{any} $\phi$ satisfying the monotonicity condition in Remark~\ref{rem:StructureForIPC}, which can be interpreted as presuming (quite plausibly) that the population's agents switch from strategy $i$ to $j$ with a rate that does not decrease when $\tilde{p}_{ij}$ increases.
\end{remark}

\subsection{Proving Theorem~1: outline and technical issues}
\label{subsec:TechnicalIssues}
In Appendix~\ref{sec:proof}, we provide a rigorous proof for Theorem~1 inspired by the approach used to establish Krasovskii-LaSalle's invariance principle~\cite{Lasalle1960Some-extensions}.

Notice that, in spite of~(\ref{eq:LeqZero}) and~(\ref{eq:DerLIneq}), we cannot employ Lyapunov's second method to prove Theorem~1 because of the following technical issues, which will also complicate our proof in Appendix~\ref{sec:proof}: {\bf (i)}~Eq.~(\ref{eq:NewEPGc}) does not directly determine an equilibrium for $q\st$. In fact, according to~(\ref{eq:NewEPGc}), $q\st$ could conceivably be unbounded. {\bf (ii)}~Challenge (i) also complicates the characterization of the equilibria of $x\st$ and its stability properties because, according to~(EDM), $q\st$ will influence the dynamics of $x\st$ via its dependence on $p\st$. {\bf (iii)}~No term accounting for the deviation $\cB\st-\beta^*$ appears on the right side of the inequality~(\ref{eq:DerLIneq}). {\bf (iv)}~Given $\check{q}$ in $\mathbb{R}$, there may be more than one best response $x$ for which $\cP(x,\check{q}\vbeta+r^o)=0$, and for a given $\check{x}$ in $\mathbb{X}$ the set $\{q \in \mathbb{R} \ | \ \cP(\check{x},q\vbeta+r^o)=0 \}$ may have more than one element and even be unbounded. 

Notice that $\cL$ would satisfy the conditions for~(\ref{eq:NewEPG}) to qualify as a $\delta$-antipassive payoff dynamics model, as defined in~\cite{Fox2013Population-Game}. However, \cite[Theorem~4.2]{Fox2013Population-Game} is not applicable to our context because of the issues (i)-(iv) listed above and the facts that the right-hand side of~(\ref{eq:DerLIneq}) is not compatible with~\cite[(41)]{Fox2013Population-Game} and we seek to characterize GAS for $\cY\st$ and not its derivative. Furthermore, because a so-called stationary game satisfying~\cite[(39)]{Park2019From-Population} clearly cannot be constructed for~(\ref{eq:NewEPG}), we cannot use the stability results in~\cite{Park2019From-Population,Park2018Payoff-Dynamic-}.

\subsection{On saturating $r\st$}

{Recall that the cost the social planner has to bear over the interval $[t,t+T]$ is equal to $\int_{t}^{t+T} r(\tau)' x(\tau) d\tau$, where $r(t) = q(t) \vec{\beta} + r^*$. Thus, when $q(t)$ is allowed to grow unbounded, so can the cost to the social planner. For this reason, if possible, it may be desirable to cap $q(t)$ until it settles near $\mathfrak{Q}^*$.}
Unfortunately, we cannot establish an a-priori bound on $q\st$, and hence on neither $r\st$ nor $p\st$. However, for certain protocols $\cT$ in~(EDM), we could replace $r\st$ with $r^{\text{\tiny SAT}}\st:= q^{\text{\tiny SAT}}\st \vbeta + r^*$, where $q^{\text{\tiny SAT}}$ would be obtained from $q$ as follows:
\begin{equation}
    q^{\text{\tiny SAT}} : = \max \{-q^{\text{\tiny MIN}},\min \{q,q^{\text{\tiny MAX}}\}\}, \quad q \in \mathbb{R}
\end{equation} Here, $q^{\text{\tiny MIN}}$ and $q^{\text{\tiny MAX}}$ would be positive constants (preferably the smallest) for which the following holds:
\begin{equation}
    \label{eq:SatEquality}
    \mathcal{V} (x, q \vbeta + r^o)=\mathcal{V} (x, q^{\text{\tiny SAT}} \vbeta + r^o), \quad q \in \mathbb{R}, \ x \in \mathbb{X}
\end{equation} For instance, for Smith's protocol discussed in Definition~\ref{def:IPC}, we could choose  $q^{\text{\tiny MAX}} =  q^{\text{\tiny MIN}} = (\bar{\mathcal{T}} / \lambda+\rho)/\min_{i \neq j} |\vbeta_i - \vbeta_j|$. Notice that~(\ref{eq:SatEquality}) would guarantee that replacing $r\st$ with $r^{\text{\tiny SAT}}\st$ would change the payoff in a way that would have no effect on $\cY\st$ and maintain the validity of Theorem~1. As $r\st'x\st$ is the normalized cost at time $t$ the social planner must bear to implement $G$ and $H$, replacing $r\st$ with $r^{\text{\tiny SAT}}\st$ would cap the cost until $q\st$ settles near~$\mathfrak{Q}^*$.

We would like to note that in our simulations in \S\ref{sec:AnytimeBounds}, in which we use Smith's protocol, we have not observed $q\st$ to be large enough to reach saturating levels.

\section{Using $\cL$ to obtain anytime bounds}
\label{sec:AnytimeBounds}

We start by using~(\ref{eq:DerLIneq}) to state~(a) below, and (b) follows from~(\ref{eq:defcL}), for $t\geq 0$:
\begin{equation}
\label{eq:LyapunovBound}
\alpha:= \cL(\cY\sta{0}) \overset{(a)}{\geq} \cL(\cY\st) \overset{(b)}{\geq} \sS(\cI\st,\cR\st,\cB\st).
\end{equation}  

Although in \S\ref{subsec:Lyapunov} we adopted the convention of not indicating in our notation that $\sS$ depends on the design parameters, here it will be useful to remember that it does and this includes dependence on $\upsilon$. In fact, we will soon outline a method for selecting $\upsilon>0$ based on~(\ref{eq:LyapunovBound}).

Assuming that $\rho^*$ and $c^*$ are pre-selected, while $\upsilon$ can vary to meet an overshoot specification, we now proceed to construct an upper-bound for $I\st/{I^*}$. Obtaining an upper bound for $I\st/{I^*}$ is important because, although Theorem~1 guarantees that $\cY\st$ will converge to $\mathbb{E}^*$, the theorem offers no guarantees on the transient behavior of $I\st$. Using $\pi^*_{\upsilon}(\alpha)$ defined below, we can leverage~(\ref{eq:LyapunovBound}) to obtain the anytime bound $I\st \leq I^* \pi^*_{\upsilon}(\alpha)$, $t \geq 0$.

\begin{definition} Given the parameters specifying~(EDM) and~(\ref{eq:NewEPG}) we seek to perform the following optimization:
\begin{equation} \label{eq:piStarDef}
\pi^*_{\upsilon}(\alpha): = \tfrac{1}{I^*} \sup \{ \ \cB^{-1}\cI  \ | \ \sS(\cY) \leq \alpha, \ \cY \in \mathbb{Y} \ \},
\end{equation} where we reverse~(\ref{eq:param-a}) to write $I=\cI / \cB$ and $I^*=\cI^* / \beta^*$. Using the fact that $\sS$ is convex~(see~Remark~\ref{rem:PropOfsS}), we conclude that~(\ref{eq:piStarDef}) is a quasi-convex program that can be swiftly solved using available software.
\end{definition}

From Remark~\ref{rem:PropOfsS}, we can immediately conclude that for any given $\upsilon>0$, it holds that $\pi^*_{\upsilon}(0)=1$ and $\pi^*_{\upsilon}(\alpha)$ is an increasing continuous function of $\alpha \geq 0$.


\subsection{Bounds: initial endemic equilibrium (n=2)}
\label{subsec:Boundsnequal2}
Throughout this subsection, consider that $n=2$ and $\cY\sta{0}$ is an endemic equilibrium point for which $(\cI,\cR)\sta{0}=(\hat{\cI},\hat{\cR})\sta{0}$, $\cB\sta{0}=\vbeta'x\sta{0}=:\beta^o$, and $q\sta{0}=0$.
 Namely, the system starts at an equilibrium that could have resulted from the prior use of~(\ref{eq:NewEPG}) in which $\beta^o$ (instead of $\beta^*$) was the endemic transmission rate. We proceed by observing that, since the entries of $p\sta{0}=r^o$ are identical (both are equal to $-c_n$), any $x(0) \in \mathbb{X}$ is in $\sM(p\sta{0})$, which implies $\cS(x(0),p\sta{0})=0$. Hence, in this case, by direct substitution into~(\ref{eq:defcL}), we obtain $\alpha=\cL(\cY\sta{0}) = \tfrac{1}{2} \upsilon^2 (\beta^o-\beta^*)^2$, which using~(\ref{eq:LyapunovBound}) leads, for $t\geq 0$,~to the following inequalities:
\begin{equation}
\label{eq:BetaAnytimeBound}
 \frac{\upsilon^2(\cB\st -\beta^*)^2}{2} \overset{(a)}{\leq} \sS(\cY\st) \overset{(b)}{\leq} \frac{\upsilon^2 \tilde{\beta}^2}{2}, \ \tilde{\beta}:=\beta^o-\beta^*
\end{equation}

Based on (a) and (b) in~(\ref{eq:BetaAnytimeBound}), it readily follows that ${| \cB\st-\beta^* | \leq | \tilde{\beta} |}$ and if $\beta^* < \beta^o$ then $\cB\st \leq \beta^o$. 

From (b) in~(\ref{eq:BetaAnytimeBound}), we also obtain: 
\begin{equation} \label{eq:AnytimeBoundnis2}
    I\st \leq I^* \times \pi^*_{\upsilon}( \tfrac{1}{2}\upsilon^2 \tilde{\beta}^2), \quad t \geq 0
\end{equation}

\begin{remark} \label{rem:AnytimeConstraint}
{\bf (Universality of~(\ref{eq:AnytimeBoundnis2}))} Analogously to Remark~\ref{rem:thm1Universality}, it is pertinent to observe that since the computation of $\pi^*_{\upsilon}(\alpha)$ for a given $\alpha$ does not require knowledge of the protocol $\cT$, (\ref{eq:AnytimeBoundnis2}) remains valid for any (EDM) satisfying the conditions of Theorem~1. Notice that obtaining a bound such as this would not have been possible in the absence of a Lyapunov function, as would have been the case for the na\"{i}ve approach described in Remark~\ref{rem:naive}.
\end{remark}

The following proposition indicates that $\upsilon$ plays a key role in bounding the overshoot of $I\st / I^*$.

\begin{proposition} \label{prop:propertiespistar}  {\bf (i)} For any $\check{\upsilon} \geq \upsilon > 0$, it holds that $\pi^*_{\check{\upsilon}}(\tfrac{1}{2} \check{\upsilon}^2 \tilde{\beta}^2 ) \geq \pi^*_{\upsilon}(\tfrac{1}{2} \upsilon^2 \tilde{\beta}^2) $. {\bf (ii)} Furthermore, it also holds that $ \pi^*_{\upsilon}(\tfrac{1}{2} \upsilon^2 \tilde{\beta}^2 ) \geq \tfrac{1}{I^*} \eta(1-\sigma\bar{\beta}^{-1}) > 1$, where we define ${\bar{\beta}:=\min \{ |\tilde{\beta}|+\beta^*,\vbeta_2 \}}$.
\end{proposition}
\noindent {\bf Proof:}  Express the constraint defining $\pi^*_{\upsilon}(\tfrac{1}{2} \upsilon^2 \tilde{\beta}^2)$ as $\sS(\cY) - \tfrac{1}{2} \upsilon^2 ( \cB-\beta^*)^2 \leq \tfrac{1}{2} \upsilon^2 (\tilde{\beta}^2 - ( \cB-\beta^*)^2)$, where for any given $\cB$ the left-hand side is a (convex) function of $(\cI,\cR)$. Using the same steps leading to (a)-(b) in~(\ref{eq:BetaAnytimeBound}) we infer that ${\tilde{\beta}^2 - ( \cB-\beta^*)^2 \geq 0}$. Hence, we can establish (i) by observing that increasing $\upsilon$ does not tighten the constraint. To show (ii) it suffices to select ${\cI=\hat{\cI}}$, ${\cR=\hat{\cR}}$ and $\cB=\bar{\beta}$ as a feasible solution. $\square$

\begin{figure}
    \centering
    \input{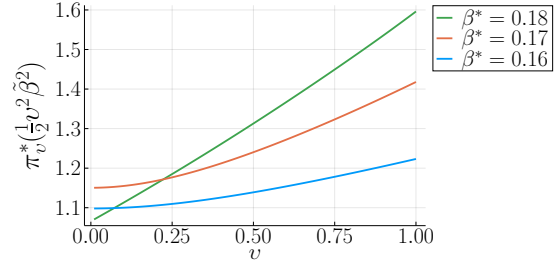} 
    \caption{Plot of $\pi^*_{\upsilon}(\tfrac{1}{2}\upsilon^2\tilde{\beta}^2)$ in Example~1 as a function of $\upsilon$ for varied $\beta^*$  (other parameters of Example~1 are unchanged).}
    \label{fig:pistar}
\end{figure}

We will use the following example to illustrate the validity of our bounds. \underline{Our time unit will be one day.}

\begin{example} \label{example1} Consider that $g=0$, $\sigma=0.1$ (infectiousness period $\sim$ 10 days), $\gamma=\sigma$, and $\omega=0.005$~(immunity period $\sim$ 200 days). The problem parameters are $\vbeta_1= 0.15$, $\vbeta_2=0.19$, while the cost vector is $c_1= 0.2$, $c_2=0$. We select $c^*=0.1$, which gives $\beta^*=0.17$, $x^*_1=x^*_2=0.5$, and $(I^*,R^*) \approx (1.96 \%,39.22 \%) $. We assume that  $x_1\sta{0}=1$, $(I\sta{0},R\sta{0}) = (\hat{I}\sta{0},\hat{R}\sta{0}) = (1.60\%,31.75 \%)$, and $\beta^o=\cB\sta{0}=0.15$. Our goal is to design $G$ and $H$ so that ${I\st} \leq 1.344 \times {I^*}$. Since $n=2$, $\rho^*$ is irrelevant (see \S\ref{subsubsec:RulesForRho}) and we can use~(\ref{eq:AnytimeBoundnis2}) to select~$\upsilon$.
\end{example}

\begin{figure}
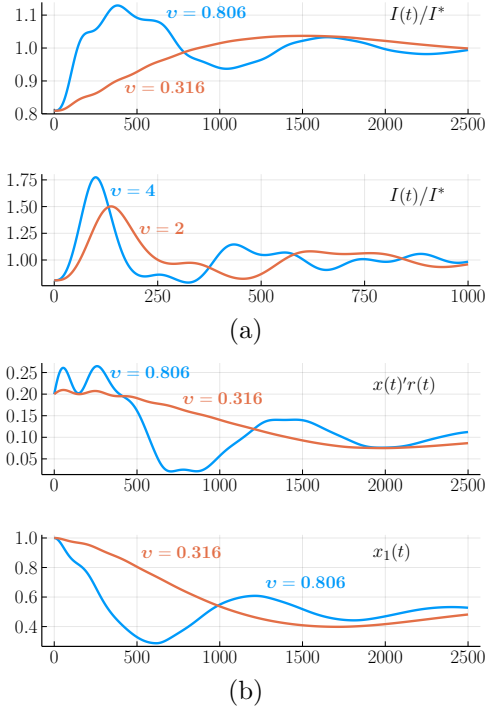

    \centering
    \input{fig1_v3.tikz} \\
    (a) \\ \vspace{0.08in}
    \input{fig2_v3.tikz} \\
    (b)
    \caption{Simulation for Example~\ref{example1} using 
    $\upsilon$ as shown, and a Smith's protocol specified by $\lambda=0.1$ and $\bar{\cT}=0.1$. }
    \label{fig:Example}
\end{figure}

Example~\ref{example1} would describe the case in which expensive measures were previously ($t<0$) in place, but a planner seeks from $t=0$ onward to relax those measures to reduce the normalized cost rate from $r'\sta{0}x\sta{0}=0.2$ to a long-term limit of $c^*=0.1$. From our numerical results (see Fig.~\ref{fig:pistar} for $\beta^*=0.17$), we determine that $\pi^*_{0.806}(\tfrac{1}{2} (0.806 \times 0.02)^2 ) \approx 1.3436$ and conclude from Proposition~\ref{prop:propertiespistar} that any positive $\upsilon \leq 0.806$ will guarantee for \emph{any} protocol satisfying the conditions of Theorem~1 that ${I\st \leq 1.344 \times I^*}$ holds, as required in Example~\ref{example1}. Fig.~\ref{fig:Example}(a) illustrates for $\upsilon \in \{0.806,0.316\}$ that the required bound indeed holds for a Smith's protocol. 
{Fig.~\ref{fig:Example}(a) suggests that (\ref{eq:AnytimeBoundnis2}) may be conservative.}
However, since~(\ref{eq:AnytimeBoundnis2}) must be valid for any protocol (not just Smith's) satisfying the conditions of Theorem~1, we do not know how conservative it may be. It is worth noting that, as illustrated in Fig.~\ref{fig:Example}(a), one could have significantly exceeded a $34.4\%$ overshoot by selecting $\upsilon \geq 2$.

From the previous discussion, one could be tempted to select a very small $\upsilon < 0.806$ expecting to perhaps eliminate any overshoot. However, as Fig.~\ref{fig:Example}(b) illustrates, smaller $\upsilon$ may lead to slower convergence, which would keep $x\st'r\st$ higher for longer. Thus, selecting the largest $\upsilon$ for which the required overshoot constraint is guaranteed by~(\ref{eq:AnytimeBoundnis2}) could be a sensible approach. In the case of Example~\ref{example1} this approach would yield $\upsilon=0.806$.




\subsection{Bounds: initial endemic equilibrium ($n\geq 3$)}
Consider that $n\geq 3$ and that $\cY\sta{0}$ satisfies the conditions specified in~\S\ref{subsec:Boundsnequal2}. Following the same argumentation as in~\S\ref{subsec:Boundsnequal2}, we conclude that $\cS(x\sta{0},p\sta{0})=0$ when the support of $x\sta{0}$ is included in that of $x^*$, in which case~(\ref{eq:BetaAnytimeBound}) remains valid. If this condition on the support of $x\sta{0}$ does not hold, then (b) in (\ref{eq:BetaAnytimeBound}) is no longer valid and we should instead use:
\begin{equation}
 \frac{\upsilon^2(\cB\st -\beta^*)^2}{2} \overset{(a)}{\leq} \sS(\cY\st) \overset{(b)}{\leq} \frac{\upsilon^2 \tilde{\beta}^2}{2} +\cS(x\sta{0},p\sta{0}),
\end{equation} where $\cS(x\sta{0},p\sta{0})$ can be computed by direct substitution when an explicit formula for $\cS$ is known. We can use~(\ref{eq:cSIPCFormula}) to compute $\cS(x\sta{0},p\sta{0})$ for any IPC protocol, such as for Smith's protocol as shown below:
\begin{align*}
    \cS^{\text{\tiny Smith}}(x(0),p(0)) = & \sum_{i=1,j=1}^{n,n} x_i\sta{0} \Phi(r^o_j-r^o_i) \\ \Phi(\nu):=&\begin{cases} \tfrac{1}{2} [\nu]_+^2 & \text{if $\nu \leq \bar{\cT}$} \\ [\nu]_+\bar{\cT} & \text{if $\nu > \bar{\cT}$} \end{cases}, \quad \nu \in \mathbb{R}
\end{align*}

\section{Conclusions}
\label{sec:conclusions}

We put forth a system theoretic methodology to model and regulate the endemic prevalence of infections 
for the case where the decisions of a population of strategically interacting agents determine the epidemic transmission rate. In our framework, the agents choose from a finite set of strategies that influence the transmission rate, and they can repeatedly revise their choices to benefit from the strategies' net rewards (payoff) resulting from incentives after deducting the strategies' intrinsic costs. An evolutionary dynamic model captures the agents' preferences by quantifying the influence of the payoff on the rate with which agents adopt or abandon the strategies. Central to our paradigm is the so-called epidemic population game (EPG) formed by two dynamically coupled components: {\bf (i)}~A modification of the standard SIRS model in which the transmission rate varies as a function of the strategies' prevalences; {\bf (ii)}~A dynamic payoff mechanism (we seek to design) that generates the strategies' incentives. 

Our main result is a dynamic payoff mechanism that is guaranteed to steer the epidemic variables (via incentives to the population) to an endemic equilibrium characterized by the lowest prevalence of infections, subject to cost constraints. Using a Lyapunov function we constructed to prove convergence, we established an (anytime) upper bound for the peak size of the population's infectious portion.

\bibliographystyle{plain}        
{{\footnotesize
\bibliography{MartinsRefs,La}  
}

\appendix

\section{Proof of Theorem~1}
\label{sec:proof}

In order to prove Theorem~1, subsequently we show that $\mathbb{E}^*$ is GAS (see~Definition~\ref{def:GAS}). 

\noindent {\bf Important notes:} {\bf (i)} We assume that $\cY\so$ is arbitrarily selected in $\mathbb{Y}$ and kept fixed throughout the rest of the proof. {\bf (ii)} We will introduce notation and definitions throughout the proof as needed. {\bf (iii)} We will make extensive use of the best response map $\sM$. We will also refer repeatedly to $\mathbb{E}^*$, $\mathfrak{Q}^*$, $\zeta_1^*$ and $\zeta_2^*$ defined in the statement of Theorem~1.

\noindent {\bf Proof structure:} The method used in~\cite{Lasalle1960Some-extensions} informed our proof strategy, which follows four main steps:  {\bf Step~1:}~We show that $\mathbb{E}^*$ satisfies (GASa); {\bf Step~2:}~We prove that $\{ q\st \ | \ t\geq 0 \}$ is bounded for any $\cY\so$ in $\mathbb{Y}$; {\bf Step~3:}~We show that $\mathbb{E}^*$ satisfies (GASb); {\bf Step~4:}~We prove that $\mathbb{E}^*$ satisfies (GASc).

\noindent {\bf On why we do not use a standard LaSalle-type theorem:} Since our proof resembles portions of a typical proof for a standard LaSalle's theorem (such as~\cite[Thm~4.4]{Khalil1995Nonlinear-syste}), it is important that we provide the following list of reasons for choosing not to use such theorems in the first place: {\bf (i)}~The state-space set $\mathbb{Y}$ is positive invariant but is not compact (is not closed). Moreover, $G$ is not defined in the closure of $\mathbb{Y}$. Hence, to employ \cite[Thm~4.4]{Khalil1995Nonlinear-syste} we would need to construct a closed positively invariant subset of $\mathbb{Y}$ for each initial condition; {\bf (ii)} Because of the technical issues listed in~\S\ref{subsec:TechnicalIssues}, in order to construct $M$, as defined in~\cite[Thm~4.4]{Khalil1995Nonlinear-syste}, we would anyway need to perform most of the analysis in steps~1-4 of our proof; {\bf (iii)} We believe that our proof method is very informative because it explains in detail why $q\st$ is bounded and why $\mathbb{E}^*$ is a positively invariant limit set.

\subsection{Step 1: showing that $\mathbb{E}^*$ satisfies (GASa)}

\begin{definition} We start with the following definitions:
\begin{equation}
    \bar{x} : = \begin{bmatrix} 0 &  \cdots & 1 \end{bmatrix}, \   \underbar{x} : = \begin{bmatrix} 1 &  \cdots & 0 \end{bmatrix}   
\end{equation}
\end{definition}

\noindent {\bf Step 1 (GASa):} Define $\mathbb{A}^*$ as follows:
\begin{equation}
\mathbb{A}^* : = \{ (x,q) \ | \ x \in \sM(q\vbeta+r^o), \vbeta'x=\beta^*\}
\end{equation} or use~(\ref{informative}),~(\ref{negdef}),~and~(NS) to rewrite $\mathbb{A}^*$ as:
 \begin{equation}
 \label{eq:EquivAstar}
\mathbb{A}^* = \{ (x,q) \in \mathbb{X}\times\mathbb{R} \  | \ \cS(x,p)=\cV(x,p)=0, \vbeta'x=\beta^*\}
\end{equation}
 Since, by Remark~\ref{rem:KKTxstar}, $\beta^*$ satisfies $\vbeta_1<\beta^*<\vbeta_n$, subsequently we will be able to show the following equality:
\begin{equation}
\label{eq:AstarCases}
 \mathbb{A}^*= \{x^*\} \times \mathfrak{Q}^*
\end{equation}

We now proceed to prove~(\ref{eq:AstarCases}) for Cases~I and~II: \\
\underline{\bf Case~I:}~Let $\{i^*,i^*+1\}$ be the support of $x^*$~(see Remark~\ref{rem:KKTxstar}). If $q=0$ then we conclude, from the fact that $\{i^*,i^*+1\}$ is the support of any $x$ in $\sM(r^o)$, that $\{x \in \sM(r^o) \ | \  \beta^* = \vec{\beta} x \}=\{x^*\}$. 
If $q>0$ and $x$ is any element in $\sM(q\vec{\beta}+r^o)$, then $x_1=\ldots=x_{i^*}=0$ and, consequently, $\vec{\beta}'x > \beta^*$. Analogously, if $q<0$ then $\vec{\beta}'x < \beta^*$ for all $x$ in $\sM(q\vec{\beta}+r^o)$. \\
\underline{\bf Case~II:} Let $i^*$ be the support of $x^*$. If $-\zeta_2^* \leq q \leq \zeta_1^*$, then $\{x \in \sM(q\vec{\beta}+r^o) \ | \  \beta^* = \vec{\beta} x \}=\{x^*\}$. If { $q> \zeta_1^*$ }, then $\sM(q\vec{\beta}+r^o) = \{ \bar{x} \}$, which would not be viable for $\mathbb{A}^*$ because $\vbeta'\bar{x}=\beta_n > \beta^*$. A similar argument shows that { $q < -\zeta_2^*$} is not viable for~$\mathbb{A}^*$.

The proof for Step~1 is concluded by using~(\ref{eq:EquivAstar}) and~(\ref{eq:AstarCases}) in conjunction with~(\ref{eq:defcL}) and Remark~\ref{rem:PropOfsS}. Or, equivalently, the following holds for $\mathbb{E}^* = (\cI^*,\cR^*) \times \mathbb{A}^* $:
\begin{equation}
\label{eq:EquivLZeroInEStar}
\cY \in \mathbb{E}^* \Leftrightarrow \cL(\cY) = 0
\end{equation}

\subsection{Additional definitions and a useful lemma}
\label{subsec:AdditionalDefsAndLemma}
Start by selecting and keeping fixed throughout the proof two positive constants $\zeta_1 > \zeta_1^*$ and $\zeta_2 > \zeta_2^*$.

We proceed with the following additional definitions:
\begin{align}
    \mathbb{P}_1 : = \Big \{ \vec{\beta}+\gamma  r^o \ \Big | \  0 \leq \gamma \leq \zeta_1^{-1} \Big \}, \\ \mathbb{P}_2 : = \Big \{ -\vec{\beta}+\gamma r^o \ \Big | \ 0 \leq \gamma \leq \zeta_2^{-1} \Big \}
\end{align} 

 \begin{remark} \label{rem:BestReponsesSigma}
    Notice that if $p$ is in $\mathbb{P}_1$ then  $\sM(p) = \{ \bar{x} \}$, and if $p$ is in $\mathbb{P}_2$ then $\sM(p) = \{ \underbar{x} \}$.
 \end{remark} 

We are now ready to state the following lemma.
\begin{lemma} \label{lem:SigmaContinuity} Given any $\epsilon > 0$, there is $\delta > 0$ such that\footnote{The lemma remains valid for any norm $\| \cdot \|$ defined in $\mathbb{R}^n$.}:
\begin{align} \label{SigmaBound1}
    \max \Big \{ \| x-\bar{x} \| \ \Big | \ x \in \mathbb{X}, \min_{p \in \mathbb{P}_1 } \cP(x,p) \leq \delta \Big \} < \epsilon \\ \label{SigmaBound2}
    \max \Big \{ \| x-\underbar{x} \| \ \Big | \ x \in \mathbb{X}, \min_{p \in \mathbb{P}_2 }\cP(x,p) \leq \delta \Big \} < \epsilon
\end{align}
\end{lemma}

\noindent {\bf Proof of Lemma~\ref{lem:SigmaContinuity}:} 
We will prove the lemma by showing that assuming it was not valid would lead to a contradiction. Hence, without loss of generality, assume that there was $\epsilon^*>0$ for which no $\delta>0$ satisfying (\ref{SigmaBound1}) existed. (The case in which no $\delta>0$ satisfying (\ref{SigmaBound2}) existed would have been analogous.) In order to reach a contradiction, we start by noticing that under the assumption the following would hold:
\begin{equation}
    \| x^{(\ell)} - \bar{x} \| \geq \epsilon^* \text{ and }  \cP(x^{(\ell)},p^{(\ell)}) \leq \tfrac{1}{\ell}, \quad \ell \geq 1
\end{equation} for a sequence $(x^{(\ell)},p^{(\ell)})$ satisfying:
\begin{align}
    x^{(\ell)} &\in \arg \max \Big \{ \| x-\bar{x} \| \ \Big | \ x \in \mathbb{X}, \min_{p \in \mathbb{P}_1 } \cP(x,p) \leq \tfrac{1}{\ell} \Big \} \\
    p^{(\ell)} &\in \arg \min_{p \in \mathbb{P}_1 } \cP(x^{(\ell)},p)
\end{align} We proceed by noting that since the sequence $(x^{(\ell)},p^{(\ell)})$ would take values in the compact set $\mathbb{X} \times \mathbb{P}_1$, it would have an accumulation point $(x^*,p^*) \in \mathbb{X} \times \mathbb{P}_1$. By continuity of $\| \cdot \|$ and $\cP$, the pair $(x^*,p^*)$ would satisfy: {\bf (i)}~${\|x^* -\bar{x} \| \geq \epsilon^*}$ and {\bf (ii)} $\cP(x^*,p^*)=0$.
However, since $p^* \in \mathbb{P}_1$, we can use Remark~\ref{rem:BestReponsesSigma}, (NS), (ii) and~(\ref{negdef}) to conclude that $x^* = \{ \bar{x} \}$, which would contradict (i).$\square$

\subsection{Step~2: proving that $\{ q\st \ | \ t\geq 0 \}$ is bounded}

 Since $\mathcal{L}$ is lower bounded, we infer from~(\ref{eq:DerLIneq}) that ${(\cI\st-\hat{\cI}\st)}$ and $(\cR\st-\hat{\cR}\st)$ are square-integrable. Also, $(\cI\st-\hat{\cI}\st)$ and $(\cR\st-\hat{\cR}\st)$ are uniformly continuous (their derivatives are bounded). Thus, we conclude using Barbalat's Lemma~(see~\cite[Theorem~1]{Farkas2016Variations-on-B}):
\begin{equation} \label{eq:F1}
\lim_{t \rightarrow \infty} (\cI\st-\hat{\cI}\st)^2 + (\cR\st-\hat{\cR}\st)^2 = 0
\end{equation} Using a similar argument and~(\ref{homogeneity}), we conclude that:
\begin{equation} \label{eq:F2}
\lim_{t \rightarrow \infty} \cP(x\st,\bar{p}\st) = 0
\end{equation} where $\bar{p}\st : = \tfrac{1}{\max\{1,|q\st|\}} p\st = \tfrac{1}{\max\{1,|q\st|\}} (q\st \vec{\beta} + r^o) $. Namely, from~(\ref{eq:DerLIneq}) we infer that $\cP(x\st,\bar{p}\st)$ is integrable. Hence, in order to use Barbalat's Lemma to prove~(\ref{eq:F2}), it suffices to establish uniform continuity of $\cP(x\st,\bar{p}\st)$ (as a function of $t$).
To do so it is helpful to use the facts that {\bf (i)}~$\hat{\cI}\st \geq \eta (\beta_1-\sigma)>0$ and {\bf (ii)}~according to~(\ref{eq:DerLIneq}) $\mathcal{L}(\cY\st)$ is bounded to infer that $\ln\cI\st$~is also bounded. By (\ref{eq:GDef}) and~(\ref{eq:NewEPGc}), this implies that $\dot{q}\st$ is bounded. Consequently, $\bar{p}\st$ is Lipschitz continuous and $\dot{x}$ is bounded because $\mathcal{V}$ in (EDM) is bounded, implying that $(x,\bar{p})\st$ is uniformly continuous. From this we conclude, by recalling that $\cP$ is Lipschitz continuous, that indeed $\cP(x\st,\bar{p}\st)$ is uniformly continuous.

Define $\xi:= \upsilon^2 \min \{\beta_n-\beta^*,\beta^*-\beta_1\}$ and select $\epsilon > 0$ such that {\bf (i)}~$\upsilon^2(\vec{\beta}'x - \beta^*) > \tfrac{2}{3} \xi$ for all $x$ in $\mathbb{X}$ satisfying $\| x - \bar{x}\| < \epsilon$ and {\bf (ii)}~$\upsilon^2(\vec{\beta}'x - \beta^*) < - \tfrac{2}{3} \xi$ for all $x$ in $\mathbb{X}$ satisfying $\| x - \underbar{x}\| < \epsilon$. From Lemma~\ref{lem:SigmaContinuity}, we know that there is $\delta>0$ such that~(\ref{SigmaBound1})-(\ref{SigmaBound2}) hold. Furthermore, from~(\ref{eq:F2}) we know that there is $\kappa$ such that, for all ${t \geq \kappa}$, we have $\cP(x\st,\bar{p}\st) \leq \delta$. Consequently, we conclude that: {\bf (a)} if $q\st \geq \zeta_1$ and ${t \geq \kappa}$, then $\bar{p}\st$ is in $\mathbb{P}_1$ and ${\min_{p \in \mathbb{P}_1 } \cP(x\st,p) \leq \cP(x\st,\bar{p}\st) \leq \delta}$ and {\bf (b)} if $q\st \leq -\zeta_2$ and $t \geq \kappa$ then $\bar{p}\st$ is in $\mathbb{P}_2$ and $\min_{p \in \mathbb{P}_2 } \cP(x\st,p) \leq \cP(x\st,\bar{p}\st) \leq \delta$. Hence, combining~(i),~(a), and~(\ref{SigmaBound1}) we arrive at~(\ref{eq:F3a}), and from~(ii),~(b), and~(\ref{SigmaBound2}) we infer~(\ref{eq:F3b}).

\begin{subequations}
\label{eq:F3}
\begin{align} \label{eq:F3a}
     q(t) \geq \zeta_1 & \implies \upsilon^2(\cB(t) - \beta^*)>  \tfrac{2}{3} \xi, \ t \geq \kappa \\ \label{eq:F3b}
     q(t) \leq - \zeta_2 & \implies  \upsilon^2( \cB(t) - \beta^*) < - \tfrac{2}{3} \xi, \ t \geq \kappa
\end{align}
\end{subequations}

From (\ref{eq:GDef}),~(\ref{eq:NewEPGc}), and~(\ref{eq:F1}), and the fact that $\hat{\cI}\st \geq \eta (\beta_1 - \sigma)>0$, we can select $\underbar{t}\geq \kappa$ satisfying:
\begin{equation}
\label{eq:F4}
    \big | \dot{q}\st+\upsilon^2(\cB\st-\beta^*) \big | < \tfrac{1}{3}\xi, \quad t \geq \underbar{t}
\end{equation} 

From~(\ref{eq:F3}), and~(\ref{eq:F4}), we can finally conclude that (a) if $q\st \geq \zeta_1$ and $t \geq \underbar{t}$, then $\dot{q} < -\tfrac{1}{3} \xi$ and (b) if $q\st \leq - \zeta_2$ and $t \geq \underbar{t}$, then $\dot{q} > \tfrac{1}{3} \xi$. Hence, we can conclude that there is $\bar{t} \geq \underbar{t}$ such that the following holds:
\begin{equation} \label{eq:F5}
 -\zeta_2 \leq q\st \leq \zeta_1, \quad t \geq \bar{t}
\end{equation} {\bf (This proves that $\{q\st \ | \ t \geq 0\}$ is bounded.)}

\subsection{Step 3: showing that $\mathbb{E}^*$ satisfies (GASb)}

Subsequently, we will show by construction the existence of an accumulation point of $\{\cY\st \ | \ t\geq 0 \}$ in $\mathbb{E}^*$.

\begin{remark} \label{rem:AccPtOfQ} Before we proceed, we observe that since $\zeta_1$ and $\zeta_2$ were any arbitrarily selected constants satisfying $\zeta_1 >\zeta_1^*$ and $\zeta_2> \zeta_2^*$, we can infer from~(\ref{eq:F5}) that any accumulation point of $\{q\st  \ | \ t\geq 0\}$ must be in $\mathfrak{Q}^*$.
\end{remark}

We start by observing that continuity of~$\dot{q}\st$ and~(\ref{eq:F5}) imply that 0 is an accumulation point of ${\{\dot{q}\st \ | \ t \geq 0 \} }$. Consequently, from~(\ref{eq:GDef}),~(\ref{eq:NewEPGc}), and~(\ref{eq:F1}), and the fact that $\hat{\cI}\st \geq \eta (\beta_1-\sigma)>0$, we conclude that $(\cI^*,\cR^*,\beta^*)$ is an accumulation point of~${\{ (\cI,\cR,\cB)\st \ | \ t \geq \bar{t} \} }$. 

Let $t^{(n)}$ be a sequence of times such that $(\cI,\cR,\cB)\sta{t^{(n)}}$ converges to $(\cI^*,\cR^*,\beta^*)$. Then, the sequence \linebreak $(\cI,\cR,\cB,x,q)\sta{t^{(n)}}$ also has an accumulation point $(\cI^*,\cR^*,\beta^*,\check{x},\check{q})$ because from~(\ref{eq:F5}) we know that, for $t \geq \bar{t}$, the pair $(x,q)\st$ takes values in a compact set $\mathbb{X} \times [-\zeta_2,\zeta_1]$. We now proceed to observe that by continuity of $\mathcal{P}$ and~(\ref{eq:F2}), it must be that $\cP \big (\check{x},\tfrac{1}{ \max \{ |\check{q}|,1 \} }(\check{q} \vec{\beta}+r^o) \big )=0$. Consequently, from~(NS),~(\ref{negdef}), the fact that $\beta^* = \vec{\beta}' \check{x}$ and Remark~\ref{rem:AccPtOfQ}, we conclude that $(\check{x},\check{q})$ must be in the set $\mathbb{A}^*$ characterized in~(\ref{eq:AstarCases}) for Cases~I and~II. \\{\bf (This concludes proof that $\mathbb{E}^*$ satisfies~(GASb).)}

\subsection{Step 4: showing that $\mathbb{E}^*$ satisfies (GASc)}

In Step~3, we constructed an accumulation point $e^*$ of $\{\cY\st \ | \ t\geq 0\}$ in $\mathbb{E}^*$. Hence, there is a sequence $t^{(n)}$  such that $\lim_{n \rightarrow \infty} \cY\sta{t^{(n)}} = e^*$. However, from~(\ref{eq:EquivLZeroInEStar}) and the continuity of $\cL$, we conclude that $\lim_{n \rightarrow \infty} \cL(\cY\sta{t^{(n)}}) = \cL(e^*)=0$. Furthermore, since~(\ref{eq:DerLIneq}) guarantees that $\cL(\cY\st)$ is non-increasing, we conclude that the following holds:
\begin{equation}
\label{eq:LimcLIsZero}
\lim_{t \rightarrow \infty} \cL(\cY\st) =0
\end{equation}
Now take any candidate accumulation point $\cY^*$ in $\mathbb{Y}$. From~(\ref{eq:LimcLIsZero}) and the continuity of $\cL$ it follows that $\cL(\cY^*)=0$, which from~(\ref{eq:EquivLZeroInEStar}) implies that $\cY^*$ must be in~$\mathbb{E}^*$. \\{\bf (This concludes proof that $\mathbb{E}^*$ satisfies~(GASc).)}

\end{document}